\documentclass[12pt]{article}

\usepackage[geometry]{ifsym}
\usepackage{verbatim}
\usepackage{fancyhdr}
\usepackage{graphicx}
\usepackage{mathrsfs}
\usepackage{amssymb}
\usepackage{cite}
\usepackage{epsfig}
\usepackage{pstricks}     
\usepackage{pst-all}      
\usepackage{caption}
\usepackage{amsfonts,amssymb, amsmath, latexsym, bm}

\newtheorem{thm}{Theorem}[section]
\newtheorem{Def}[thm]{Definition}

\newtheorem{pps}[thm]{Proposition}
\newtheorem{cor}[thm]{Corollary}
\newtheorem{lem}[thm]{Lemma}

\newtheorem{pro}[thm]{Problem}
\newenvironment{pf}[1][Proof]{\noindent\textbf{#1.} }{\hfill\rule{1mm}{2mm}}

\makeatletter \@addtoreset{equation}{section} \makeatother

\begin{document}

\title{Graham's pebbling conjecture on Cartesian product of the middle graphs of even cycles\thanks{Supported by ``the Fundamental Research Funds for the Central Universities" and the NSF of the People's Republic of China(Grant
 No. 61272008, No. 11271348 and
 No. 10871189).}
  }
\author
{Zheng-Jiang Xia,\;\; Yong-Liang Pan\footnote{Corresponding
author:
ylpan@ustc.edu.cn},\quad Jun-Ming Xu,\quad Xi-Ming Cheng\\ \\
{\small     School of Mathematical Sciences,}\\
{\small             University of Science and Technology of China,}   \\
{\small             Hefei, Anhui, 230026, P. R. China}\\
{\small             Email: ylpan@ustc.edu.cn}  \\
}

\date{}
\maketitle

\date{}
\noindent\textbf{Abstract}: A pebbling move on a graph $G$ consists of taking two pebbles off one vertex and placing one on an adjacent
vertex. The pebbling number of a graph $G$, denoted by $f(G)$, is the least integer $n$ such that, however $n$
pebbles are located on the vertices of $G$, we can move one pebble to any vertex by a sequence
of pebbling moves. Let $M(G)$ be the middle graph of $G$.
For any connected graphs $G$ and $H$, Graham conjectured that $f(G\times H)\leq f(G)f(H)$.
In this paper, we give the pebbling number of some graphs
and prove that Graham's conjecture holds for the middle graphs of some even cycles.\\



\section{Introduction}

Pebbling in graphs was first introduced by Chung \cite{c89}. Consider a connected graph with a fixed number of pebbles
distributed on its vertices. A pebbling move consists of the removal of two pebbles from a vertex and the placement of
one pebble on an adjacent vertex. The pebbling number of a vertex $v$, the target vertex, in a graph $G$ is the smallest number $f(G,v)$
with the property that, from every placement of $f(G,v)$ pebbles on $G$, it is possible to move one pebble to $v$ by a sequence of
pebbling moves. The pebbling number of a graph $G$, denoted by $f(G)$, is the maximum of $f(G,v)$ over all the vertices of $G$.\\

There are some known results regarding the pebbling number (see \cite{c89,yzz12,fk01,lqwm06,sf00}). If one pebble is placed on each vertex other than the vertex $v$,
then no pebble can be moved to $v$. Also, if $u$ is at a distance $d$ from $v$, and $2^d-1$ pebbles are placed on $u$, then no pebble can
be moved to $v$. So it is clear that $f (G) \geq \max \{|V(G)|, 2^D\}$, where $D$ is the diameter of graph $G$. Furthermore, we know that
$f (K_n) = n$ and $f (P_n) = 2^n-1$ (see \cite{c89}), where $K_n$ is the complete graph  and $P_n$ is the path, respectively  on  $n$ vertices.\\

The {\it middle graph} of a graph $G$, denoted by $M(G)$, is obtained from $G$ by inserting a new vertex into each edge of $G$,
and joining the new vertices by an edge if the two edges they inserted share the same vertex of $G$.\\

\indent Given two disjoint graphs $G_1=(V_1,E_1)$ and
$G_2=(V_2,E_2)$, the Cartesian product  of them is denoted by
$G_1\times G_2$. It has vertex set $V_1\times V_2=\{(u_i,v_j)|
u_i\in V_1, v_j\in V_2\}$, where $(u_1,v_1)$ is adjacent to
$(u_2,v_2)$ if and only if $u_1=u_2$ and $(v_1,v_2)\in E_2$, or
$(u_1,u_2)\in E_1$ and $v_1=v_2$. One may view $G_1\times G_2$ as
the graph obtained from $G_2$ by replacing each of its  vertices
with a copy of $G_1$, and each of its edges with $|V_1|$ edges
joining corresponding vertices
of $G_1$ in the two copies.  Let $u\in G, v\in H$, then $u(H)$ and $v(G)$ are subgraphs of $G\times H$ with
$V(u(H))=\{(u,v)|v\in V(H)\}$, $E(u(H))=\{(u,v)(u,v')|vv'\in E(H)\}$ and $V(v(G))=\{(u,v)|u\in V(G)\}$, $E(v(G))=\{(u,v)(u',v)|uu'\in E(G)\}$.
It is clear that $u(H)\cong H$ and $v(G)\cong G$.\\

The following conjecture (see [2]), by Ronald Graham, suggests a constraint on the pebbling number
of the product of two graphs.\\

\indent Conjecture (Graham). The pebbling number of $G\times H$ satisfies $f(G\times H)\leq f(G)f(H)$.\\

Ye {\it et al.} (see \cite{yzz}) proved that  $f(M(C_{2n+1})\times M(C_{2m+1}))\leq f(M(C_{2n+1}))f(M(C_{2m+1}))$ and $f(M(C_{2n})\times M(C_{2m+1}))\leq f(M(C_{2n}))f(M(C_{2m+1})).$
In this paper, we will prove that $f(M(C_{2n})\times M(C_{2m}))\leq f(M(C_{2n}))f(M(C_{2m}))$ for $m,n\geq5$ and $|n-m|\geq2$.\\

Throughout this paper, $G$ will denote a simple connected graph with vertex set $V(G)$ and edge set $E(G)$.
$P_n$ and $C_n$ will denote a path and a cycle with $n$ vertices, respectively.
Given a distribution of pebbles on the vertices of $G$, define $p(K)$ to be the number of pebbles on a subgraph $K$ of G and $p(v)$ to be the number of pebbles on a
vertex $v$ of $G$. Moreover, we let $\tilde{p}(K)$ and $\tilde{p}(v)$ denote the numbers of pebbles on $K$ and $v$ after some sequence of pebbling moves,
respectively.


\section{Main results}

\begin{Def}{\rm(see \cite{sf00})}
Let $P_{n}=v_1v_2\cdots v_{n}$ be a path. We say that $P_n$ has weight
$\sum\limits_{i=1}^{n-1}2^{i-1}p(v_i)$ with respect to $v_n$ and this is written as $\omega_{P_n}(v_n)$.
\end{Def}

\begin{pps}{\rm{(see \cite{sf00})}}\label{pps2.1}
Let  $P_{n}=v_1v_2\cdots v_{n}$ be a path. If $\omega_{P_n}(v_n)\geq k2^{n-1}$, then at least $k$
pebbles can be moved from $P_n\backslash v_n$ to $v_n$.
\end{pps}

\begin{cor}\label{cor4.1}
Let $P_{n}=v_1v_2\cdots v_{n}$ be a path. Let $\omega_{P_n}(v_k)=\sum\limits_{i=1}^{k-1}2^{i-1}p(v_i)
+\sum\limits_{j=k+1}^{n}2^{n-j}p(v_j)$ for $2\leq k\leq n-1$.
If $\omega_{P_n}(v_k)\geq t2^{k-1}+2^{n-k}-1$ for $\frac{n+1}{2}\leq k\leq n$,
 $\omega_{P_n}(v_k)\geq2^{k-1}+t2^{n-k}-1$ for $1\leq k<\frac{n+1}{2}$, then at least $t$ pebbles can be moved from $P_{n}\backslash v_k$ to $v_k$.
\end{cor}

\begin{pf}
Without loss of generality, we assume that $\frac{n+1}{2}\leq k\leq n$.

If $k=n$, it follows from Proposition~\ref{pps2.1}.

If $\frac{n+1}{2}\leq k\leq n-1$, let $L_1=v_1v_2\cdots v_k$, $L_2=v_kv_{k+1}\cdots v_n$ be two subpaths of $P_n$.

Suppose $\omega_{P_n}(v_k)\geq t2^{k-1}+2^{n-k}-1$,
then either $\sum\limits_{i=1}^{k-1}2^{i-1}p(v_i)\geq t2^{k-1}$ or $\sum\limits_{j=k+1}^{n}2^{n-j}p(v_j)\geq2^{n-k}$ holds.

Case $1$. $\sum\limits_{i=1}^{k-1}2^{i-1}p(v_i)\geq t2^{k-1}$, by Proposition~\ref{pps2.1},
we can move $t$ pebbles from $L_1\backslash v_k$ to $v_k$.

Case $2$. $\sum\limits_{j=k+1}^{n}2^{n-j}p(v_j)\geq2^{n-k}$,
we may assume that $\sum\limits_{j=k+1}^{n}2^{n-j}p(v_j)=s2^{n-k}+h$,
where $s$ and $h$ are integers satisfying $s\geq1$ and $0\leq h< 2^{n-k}$. With $p(v_j)$ pebbles on $v_j$ $(k+1\leq j\leq n)$,
we can move $s$ pebbles from $L_2\backslash v_k$ to $v_k$.

Note that $2^{k-1}\geq2^{n-k}$ for $k\geq\frac{n+1}{2}$, we have
\begin{align*}
\sum\limits_{i=1}^{k-1}2^{i-1}p(v_i)=&\omega_{P_n}(v_k)-\sum\limits_{j=k+1}^{n}2^{n-j}p(v_j)\\
\geq&t2^{k-1}+2^{n-k}-1-(s2^{n-k}+h)\\
=&(t2^{k-1}-s2^{n-k})+(2^{n-k}-h)-1\\
\geq&(t-s)2^{k-1}.
\end{align*}

So we can move $t-s$ pebbles from $L_1\backslash v_k$ to $v_k$ with $p(v_i)$ pebbles on $v_i$ $(1\leq i\leq k-1)$.
That is to say we can move $s+(t-s)=t$ pebbles to $v_k$.
\end{pf}

\begin{cor}\label{cor2.1}
Let $P_n=v_1v_2\cdots v_n$ be a path. Then $f(M(P_n)-\{v_1,v_n\})=2^{n-2}+n-2$.
\end{cor}

\begin{figure}[ht]
\begin{center}
\hspace*{30pt}
\begin{pspicture}(-5,0)(5,3)
\cnode(-3.3,2){3pt}{u1}\cnode(-2,2){3pt}{u2}\cnode(-1,2){3pt}{u3}\cnode(0,2){3pt}{u4}\cnode(1,2){3pt}{u5}\cnode(2,2){3pt}{u6}\cnode(3.3,2){3pt}{u7}
\cnode(-1.5,1){3pt}{v1}\cnode(-0.5,1){3pt}{v2}\cnode(1.5,1){3pt}{v3}
\ncline[linestyle=dashed,dash=3pt 2pt]{u1}{u2}\ncline{u2}{u3}\ncline{u3}{u4}\ncline{u5}{u6}\ncline[linestyle=dashed,dash=3pt 2pt]{u6}{u7}
\ncline{u2}{v1}\ncline{u3}{v1}\ncline{u3}{v2}\ncline{u4}{v2}\ncline{u5}{v3}\ncline{u6}{v3}
\rput(.5,2){$\bm{\cdots}$}

\rput(-3,2.35){$v_{_1}$}\rput(-2,2.35){$u_{_1}$}\rput(-1,2.35){$u_{_{2}}$}\rput(0,2.35){$u_{_{3}}$}\rput(1,2.35){$u_{_{n-2}}$}
\rput(2,2.35){$u_{_{n-1}}$}\rput(3,2.35){$v_{_{n}}$}
\rput(-1.5,0.65){$v_{_2}$}\rput(-.5,0.65){$v_{_3}$}\rput(1.5,0.65){$v_{_{n-1}}$}
\psellipse(0,1.7)(2.7,1.65)

%
\end{pspicture}
\caption{\small The graph $M(P_n)-\{v_1,v_n\}$ in Corollary~\ref{cor2.1}. \label{fig1}
}
\end{center}
\end{figure}
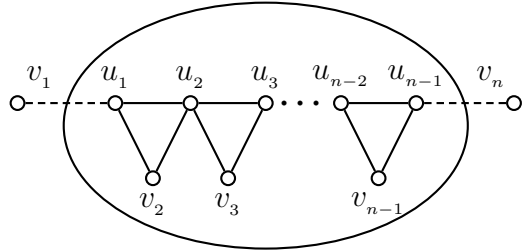

\begin{pf}
To get $M(P_n)$, we insert $u_i$ into the edge $v_iv_{i+1}$ and add the edge $u_iu_{i+1}$ for each $i\in\{1,2,\ldots ,n-2\}$.
Let $U=u_1u_2\cdots u_{n-1}$ be a subpath of $M(P_n)-\{v_1,v_n\}$.

It is clear that $f(M(P_n)-\{v_1,v_n\})\geq2^{n-2}+n-2$.
If we place one pebble on each of vertices $v_2,\ldots ,v_{n-1}$, and place $2^{n-2}-1$ pebbles on $u_{n-1}$, then we can not move one pebble to $u_1$.
So $f(M(P_n)-\{v_1,v_n\})\geq2^{n-2}+n-2$.

Now, assume that $2^{n-2}+n-2$ pebbles are located at $V(M(P_n)-\{v_1,v_n\})$.

First, we prove that one pebble can be moved to $u_k$ $(1\leq k\leq n-1)$.

While $m\leq k$, we can move $\lfloor p(v_m)/2\rfloor$ pebbles from $v_m$ to $u_m$.
While $m> k$, we can move $\lfloor p(v_m)/2\rfloor$ pebbles from $v_m$ to $u_{m-1}$.
\begin{align*}
\omega_U(u_k)\geq& 2^{n-2}+n-2-\sum\limits_{t=2}^{n-1}p(v_t)+2\sum\limits_{t=2}^{n-1}\lfloor p(v_t)/2\rfloor\\
\geq&2^{n-2}.
\end{align*}

It is clear that $2^{n-2}\geq 2^{k-1}+2^{n-k-1}-1$ for $1\leq k\leq n-1$.
By Corollary~\ref{cor4.1}, we can move one pebble from $U\backslash u_k$ to $u_k$ $(1\leq k\leq n-1)$.

Now we prove that one pebble can be moved to $v_k$ $(2\leq k\leq n-1)$.
Without loss of generality, we assume that $k\geq\frac{n+1}{2}$.

While $m< k$, we can move $\lfloor p(v_m)/2\rfloor$ pebbles from $v_m$ to $u_m$.
While $m> k$, we can move $\lfloor p(v_m)/2\rfloor$ pebbles from $v_m$ to $u_{m-1}$.

We will prove that after a sequence of pebbling moves above, two pebbles can be moved from $U$ to $u_{k-1}$,
so that one pebble can be moved from $u_{k-1}$ to $v_k$.

We consider the worst case, that is $p(u_{k-1})=0$.
\begin{align*}
\omega_U(u_{k-1})\geq&2^{n-2}+n-2-\sum\limits_{j=2\atop j\neq k}^{n-1}p(v_j)+2\sum\limits _{j=2\atop j\neq k}^{n-1}\lfloor p(v_j)/2\rfloor\\
\geq&2^{n-2}+1.
\end{align*}
It is clear that $2^{n-2}+1\geq 2\times2^{(k-1)-1}+2^{n-(k-1)-1}-1$ for $\frac{n-1}{2}\leq k-1\leq n-2$.
By Corollary~\ref{cor4.1}, we can move two pebbles from $U\backslash u_{k-1}$ to $u_{k-1}$ $(\frac{n-1}{2}\leq k-1\leq n-2)$.
So we can move one pebble to $v_k$ $(\frac{n+1}{2}\leq k\leq n-1)$, and we are done.
\end{pf}

\begin{Def}{\rm{(see \cite{sf00})}}
The $t$-pebbling number of a graph $G$ is the smallest number $f_t(G)$
with the property that from every placement of $f_t(G)$ pebbles on $G$, it is possible to
move $t$ pebbles to any vertex $v$ by a sequence of pebbling moves.
\end{Def}

\begin{lem}{\rm(see \cite{yzz})}\label{lem2.3}
If $n\geq 2$, then $f(M(C_{2n}))=2^{n+1}+2n-2$.
\end{lem}

\begin{cor}\label{cor2.2}
If $n\geq 2$, then $f_t(M(C_{2n}))\leq t2^{n+1}+2n-2$.
\end{cor}

\begin{pf}
Let $C_{2n}=v_0v_1\cdots v_{2n-1}v_0$, $M(C_{2n})$ is obtained from $C_{2n}$ by inserting $u_i$ into $v_iv_{(i+1)mod(2n)}$,
and connecting $u_iu_{(i+1)mod(2n)}$ $(0\leq i \leq2n-1)$.

Without loss of generality, we may assume that our target vertex is $u_0$ or $v_0$.

Case $1$. The target vertex is $u_0$. In this case, we use induction on $t$.

The result is obvious for $t=1$  from Lemma~\ref{lem2.3}.

Now suppose that $t2^{n+1}+2n-2$ pebbles are located at the vertices of $M(C_{2n})$.

We consider the worst case, that is $p(u_0)=0$.

Let $A=\{u_0,v_1,u_1,\ldots ,v_n,u_n\}$,
$B=\{u_n,v_{n+1},\ldots ,v_{2n-1},u_{2n-1},v_0,u_0\}$ and $G= M(C_{2n})$.
 Then we have either $A$ or $B$ contains more than $2^n+n$ pebbles.

 Note that $G[A]\cong G[B]\cong M(P_{n+2})-\{v_1,v_{n+2}\}$, according to Corollary~\ref{cor2.1},
 with $2^n+n$ pebbles on $A$ or $B$, one pebble can be moved to $u_0$.

 Note that $2^n+n\leq2^{n+1}$,
 the number of remaining pebbles is more than $(t-1)2^{n+1}+2n-2$.
 So we can move $t-1$ pebbles to $u_0$ with the remaining pebbles by induction, and we are done.

Case $2$. The target vertex is $v_0$.

Let $A'=\{u_0,v_1,\ldots ,v_{n-1},u_{n-1}\}$, $B'=\{u_{2n-1},v_{2n-1},\ldots ,v_{n+1},u_n\}$.

Suppose that $t2^{n+1}+2n-2$ pebbles are located at the vertices of $M(C_{2n})$.

We consider the worst case, that is $p(v_0)=0$.

By proposition~\ref{pps2.1}, while $p(v_n)\geq t2^{n+1}$, $t$ pebbles can be moved to $v_0$.

Now suppose that $t2^{n+1}-h$ pebbles are located at $v_n$,
without loss of generality, we assume that $p(A')\geq p(B')$, that is $p(A')\geq n-1+\lceil h/2\rceil$.

Let $L=v_0u_0u_1\cdots u_{n-1}v_n$ be a subpath of $G$ with length $n+1$ and $q=\sum\limits_{i=0}^{n-1}p(u_i)$.

While $q\geq\lceil h/2\rceil$,

$$\omega_L(v_0)=p(v_n)+\sum_{i=0}^{n-1}2^{n-i}p(u_i)\geq t2^{n+1}-h+2q\geq t2^{n+1}.$$

By Proposition~\ref{pps2.1}, $t$ pebbles can be moved from $L\backslash v_0$ to $v_0$.

While $q<\lceil h/2\rceil$, then $\sum\limits_{j=1}^{n-1}p(v_j)\geq n-1+\lceil h/2\rceil-q$. So we can move at least $\left\lfloor\frac{1}{2}(\lceil\frac{h}{2}\rceil+1-q)\right\rfloor$ pebbles to the set $\{u_0,u_1,\ldots ,u_{n-2}\}$.
Then we have
$$\omega_L(v_0)=p(v_n)+\sum_{i=0}^{n-1}2^{n-i}\tilde{p}(u_i)\geq t2^{n+1}-h+2q+4\times\frac{1}{2}(\frac{h}{2}-q)\geq t2^{n+1}.$$
By Proposition~\ref{pps2.1}, $t$ pebbles can be moved from $L\backslash v_0$ to $v_0$. The result follows.
\end{pf}

\begin{thm}\label{thm2.2}
If $m,n\geq5$ and $|n-m|\geq2$, then
$$f(M(C_{2n})\times M(C_{2m}))\leq f(M(C_{2n}))f(M(C_{2m})).$$
\end{thm}

\begin{pf}
Without loss of generality, we assume that $n\geq m+2$ $(m\geq5)$.

Let $V(M(C_{2n}))=\{u_1,u_2,\ldots ,u_{4n}\}$, $V(M(C_{2m}))=\{v_1,v_2,\ldots ,v_{4m}\}$.

For simplicity, let $G= M(C_{2n})\times M(C_{2m})$.

Now assume $(2^{n+1}+2n-2)(2^{m+1}+2m-2)$ pebbles have been placed arbitrarily at the vertices of $G$.

We may assume our target vertex is $(u_i,v_j)$, then $(u_i,v_j)$ belongs to both $u_i(M(C_{2m}))$ and $v_j(M(C_{2n}))$.

If $p(u_i(M(C_{2m})))\geq2^{m+1}+2m-2$ or $p(v_j(M(C_{2n})))\geq2^{n+1}+2n-2$, we can move one pebble to $(u_i,v_j)$ by lemma~\ref{lem2.3}.

Suppose that $p(u_i(M(C_{2m})))\leq2^{m+1}+2m-3$ and $p(v_j(M(C_{2n})))\leq2^{n+1}+2n-3$.

We will prove that if we move as many as possible pebbles from $u_l(M(C_{2m}))$ to $(u_l,v_j)$
which belongs to $v_j(M(C_{2n}))$ $(1\leq l\leq 4n)$,
 then one pebble can be moved from $v_j(M(C_{2n}))$ to $(u_i,v_j)$.

We may assume that
$$p_k=p(u_k(M(C_{2m})))\leq2^{m+1}+2m-3~(1\leq k\leq s)$$
and
$$p_k=p(u_k(M(C_{2m})))\geq2^{m+1}+2m-2~(s+1\leq k\leq 4n).$$

Now we consider the worst case scenario (i.e. the most wasteful distribution of pebbles possible).
Therefore we may assume that
\begin{align*}
p_k=\left\{
\begin{array}{ll}
2^{m+1}+2m-3 & 1\leq k\leq s,\\
t_k2^{m+1}+2m-2+(2^{m+1}-1) & s+1\leq k\leq 4n-1,\\
t_k2^{m+1}+2m-2+R & k=4n,
\end{array}
\right.
\end{align*}
where $0\leq R\leq2^{m+1}-1$ and $t_k$ is a positive integer.
According to Corollary~\ref{cor2.2}, we can move at least $\sum\limits_{k=s+1}^{4n}t_k$ pebbles to $v_j(M(C_{2n}))$.

Let
\begin{align*}
\Delta =&~(2^{n+1}+2n-2)(2^{m+1}+2m-2)-s(2^{m+1}+2m-3)\\
      &~-(4n-s-1)(2^{m+1}-1)-(4n-s)(2m-2)\\
= &~(2^{n+1}-2n-2)(2^{m+1}+2m-2)+2^{m+1}+4n-1.
\end{align*}

Therefore,
 $$\frac{\Delta}{2^{m+1}}=2^{n+1}-2n-1+\frac{1}{2^{m+1}}\left[(2^{n+1}-2n-2)(2m-2)+4n-1\right].$$

Note that $\Delta=\left(\sum\limits_{k=s+1}^{4n}t_k\right)2^{m+1}+R$, so $\sum\limits_{k=s+1}^{4n}t_k>\frac{\Delta}{2^{m+1}}-1$.
It follows that
\begin{align*}
p(v_j(M(C_{2n})))\geq~\sum_{k=s+1}^{4n}t_k
      >~2^{n+1}-2n-2+\frac{1}{2^{m+1}}\left[(2^{n+1}-2n-2)(2m-2)+4n-1\right].
\end{align*}

To the end, we only need to prove that we can move one pebble from
$v_j(M(C_{2n}))$ to $(u_i,v_j)$ with $2^{n+1}-2n-2+\frac{1}{2^{m+1}}\left[(2^{n+1}-2n-2)(2m-2)+4n-1\right]$ pebbles.

So we only need to prove that
$$2^{n+1}-2n-2+\frac{1}{2^{m+1}}\left[(2^{n+1}-2n-2)(2m-2)+4n-1\right]\geq2^{n+1}+2n-2,$$
that is
\begin{equation}\label{eq1}
2^{m+1}<\frac{m-1}{n}(2^n-1)-m+2.
\end{equation}

For $n\geq m+2\geq7$, it is clear that the right side of (\ref{eq1}) is an increasing function of $n$.
So we only need to prove (\ref{eq1}) under $n=m+2$. Substituting $n=m+2$ into (\ref{eq1}), we have
$$2^{m+1}<\frac{m-1}{m+2}(2^{m+2}-1)-m+2,$$
that is
\begin{equation}\label{eq2}
(2m-8)2^m-m^2-m+5>0.
\end{equation}

The left side of (\ref{eq2}) is an increasing function of $m$ while $m\geq5$.
When $m=5$, (\ref{eq2}) holds. This completes the proof.
\end{pf}


\section{Remark}

In fact, by a similar processing as in the proof of Corollary~\ref{cor2.2}, for any $u\in M(C_{2n})$ but $u\not\in C_{2n}$, we can prove that
\begin{cor}
If $n\geq 2$, then $f_t(M(C_{2n}),u)\leq2^{n+1}+2n-2+(t-1)(2^n+n).$
\end{cor}

Then we can prove the following theorem.

\begin{thm}
If $(u,v)\not\in C_{2n}\times C_{2m}$, where $C_{2n}\times C_{2m}$ is a subgraph of $M(C_{2n})\times M(C_{2m})$, then
$$f(M(C_{2n})\times M(C_{2m}),(u,v))\leq f(M(C_{2n}))f(M(C_{2m})).$$
\end{thm}

\begin{pf}
If $(u,v)\not\in C_{2n}\times C_{2m}$, then we can get $u(M(C_{2m}))\nsubseteq C_{2n}\times M(C_{2m})$ or
$v(M(C_{2n}))\nsubseteq M(C_{2n})\times C_{2m}$.

Without loss of generality, we assume that $u(M(C_{2m}))\nsubseteq C_{2n}\times M(C_{2m})$.

Let $V(M(C_{2m}))=\{v_1,v_2,\ldots,v_{4m}\}$.

If we move as many as possible pebbles from $v_j(M(C_{2n}))$ to $(u,v_j)$ which belongs to $u(M(C_{2m}))$ $(1\leq j\leq 4m)$,
by a similar processing as in the proof of Theorem~\ref{thm2.2}, we can prove that the number of pebbles on $u(M(C_{2m}))$
is more than $2^{m+1}+2m-2$, so one pebble can be moved from $u(M(C_{2m}))$ to $(u,v)$ with these pebbles.
\end{pf}

In this paper, we have shown that while $m,n\geq5$ and $|m-n|\geq2$, $f(M(C_{2n})\times M(C_{2m}))\leq f(M(C_{2n}))f(M(C_{2m}))$.
The remaining question is  open.
\begin{pro}
$f(M(C_{2n})\times M(C_{2m}))\leq f(M(C_{2n}))f(M(C_{2m}))$, for $m=n$ or $m=n-1$.
\end{pro}

\end{document}